\documentclass[11pt,leqno]{article}
\usepackage{graphicx, amsfonts, amsthm, amsxtra, amssymb, verbatim, makeidx}
\usepackage{subeqnarray, relsize}
\usepackage[mathscr]{euscript}
\usepackage{hyperref}
\usepackage{color} %
\hypersetup{
    colorlinks=true,       
    linkcolor=blue,          
    citecolor=blue,        
    filecolor=blue,      
    urlcolor=blue           
}\usepackage[nameinlink,capitalize,noabbrev]{cleveref}
\makeindex
\makeglossary
\usepackage{amsmath,amssymb, amsthm}
\usepackage[utf8]{inputenc}
\usepackage[activeacute, english]{babel}
\usepackage{wrapfig}
\usepackage{amssymb, amsmath, amsthm} 
\usepackage{graphicx}
\usepackage{color}
\usepackage{hyperref}
\usepackage{amssymb}
\usepackage{url}
\usepackage{pdfpages}
\usepackage{fancyhdr}
\usepackage{hyperref}
\usepackage{subfig}
\usepackage{mathrsfs}
\usepackage{float}
\usepackage{pgf,tikz}
\usetikzlibrary{arrows}
\usetikzlibrary[patterns]
\numberwithin{equation}{section}

\DeclareMathOperator{\pn}{{\sf p}}

\newcommand{\N}{\mathbb{N}}
\newcommand{\Z}{\mathbb{Z}}
\newcommand{\C}{\mathbb{C}}
\newcommand{\Q}{\mathbb{Q}}

\renewcommand{\P}{\mathbb{P}}

\newcommand{\fB}{\mathfrak{B}}
\newcommand{\fD}{\mathfrak{D}}

\newcommand{\fC}{\mathfrak{C}}
\newcommand{\cE}{\mathcal{E}}

\newcommand{\cK}{\mathcal{K}}

\newcommand{\bone}{\mathbf{1}}

\newcommand{\res}{\text{res}}

\newcommand{\prim}{\text{prim}}

\newtheorem{thm}{Theorem}[section]

\newtheorem{conj}{Conjecture}

\newtheorem{lemma}[thm]{Lemma}

\newtheorem{cor}[thm]{Corollary}

\newtheorem{prop}[thm]{Proposition}

\theoremstyle{definition}

\newtheorem{rmk}{Remark}[section]

\newtheorem{dfn}{Definition}[section]

\renewenvironment{proof}{{\bfseries \noindent Proof} }{ \qed \\}

\newcount\colveccount 
\newcommand*\colvec[1]{ 
    \global\colveccount#1
    \begin{pmatrix} 
    \colvecnext 
}

\def\colvecnext#1{
    #1
    \global\advance\colveccount-1 
    \ifnum\colveccount>0 
        \\
        \expandafter\colvecnext 
    \else 
        \end{pmatrix}
    \fi
}

\usepackage{tikz-cd}
\usepackage{booktabs}

\begin{document}

\begin{center}
{\LARGE\bf  Hodge cycles of Klein hypersurfaces}\\

{{\sc Mario Pastrana and Roberto Villaflor}}

\end{center}
 
\begin{abstract}
We study the space of Hodge cycles on Klein varieties. Relating it to the Hodge characters in a covering Fermat variety, we decompose the space of Hodge characters into a part generated by classes of rational subvarieties, and a non-trivial part which we call exceptional. Using Dirichlet characters we show that the exceptional part is always finite.
Furthermore, we give sufficient conditions under which there are no exceptional characters, providing infinitely many new cases of the Hodge conjecture for Klein varieties. Finally, we describe the set of exceptional characters of Klein varieties for
small dimensions.
\end{abstract}

\def\pn{{\sf p}}
\def\rootsG{{\sf G}}
\def\NLL{{\rm NL}}   

\section*{Introduction}
\label{intro}
Fermat varieties $X^n_m=\{x_0^m+\cdots+x_{n+1}^m=0\}\subseteq\P^{n+1}$ constitute the most studied case of the Hodge conjecture \cite{Shi79,Ran1980,Shioda1981,AokiShioda1983,Aoki1983,aoki1987,Schuett2010,kang2016refined}, and still remains widely open (see for instance \cite{MMRV2026}). These varieties have a big abelian automorphism subgroup $G^n_m=\mu_m^{n+2}/\Delta(\mu_n)$ acting diagonally on $X^n_m$ by multiplication of $m$-roots of unity. This action induces a spectral decomposition on the primitive cohomology, which refines the Hodge decomposition \cite[Theorem 1]{Shi79} in terms of one dimensional eigenspaces associated to its characters, which we denote $$\widehat{G}^n_m=\left\{\alpha=(a_0,\ldots,a_{n+1})\in(\Z/m)^{n+2}: \sum_{i}a_i=0\right\}.$$
More precisely, letting $\mathfrak{A}^n_m=\{\alpha\in\widehat{G}^n_m: a_i\neq 0\ \forall i\}$ we have that $$\mathrm{H}^n(X^n_m,\C)_\prim=\bigoplus_{\alpha\in\mathfrak{A}^n_m}V(\alpha)$$ and, defining $|\alpha|:=\frac{1}{m}\sum_{i}\overline{a_i}\in\{1,\ldots,n+1\}$ (where $\overline{a_i}\in\{0,\ldots,m-1\}$ denotes its remainder), we have for $p+q=n$ that 
$$
\mathrm{H}^{p,q}(X^n_m)_\prim=\bigoplus_{|\alpha|=q+1}V(\alpha).
$$
Since the spectral decomposition is defined over the cyclotomic field $\Q(\zeta_m)$, its Galois group $\text{Gal}(\Q(\zeta_m)/\Q)\simeq(\Z/m)^\times$ acts permuting the eigenspaces and so, when $n=2p$ is even, we can describe the primitive Hodge cycles as
$$
\mathrm{Hod}^n(X^n_m)_\prim\otimes\C=\bigoplus_{\alpha\in\mathfrak{B}^n_m}V(\alpha),
$$
where $\mathfrak{B}^n_m=\{\alpha\in\mathfrak{A}^n_m:|t\cdot \alpha|=p+1 \ \forall t\in(\Z/m)^\times\}$. This allows one to reduce the study of Hodge cycles to understand $\mathfrak{B}^n_m$, which amounts to study a set of characters with elementary arithmetic properties. Furthermore, if we denote by $\mathfrak{C}^n_m\subseteq\fB^n_m$ the set of characters supporting algebraic classes, one can reduce the Hodge conjecture to verify when $\mathfrak{C}^n_m=\fB^n_m$ (see \cite[Proposition 3.1]{MMRV2026}). The set of characters supported in linear cycles is denoted by $\mathfrak{D}^n_m\subseteq\mathfrak{C}^n_m$, and the theorem of Aoki \cite[Theorem A]{Aoki1983} determines precisely when they generate all Hodge characters
\[
\fD_m^n=\fB_m^n
\quad\Longleftrightarrow\quad
m=4,\ \text{ or } m \text{ prime, or } \gcd\big(m,(n+1)!\big)=1 ,
\]
providing infinitely many cases where $X^n_m$ satisfies the Hodge conjecture.

On one hand the spectral refinement is the key advantage to study the Hodge conjecture on $X^n_m$, but at the same time its high symmetry also provides a lot of algebraic cycles inside them, and the challenge for general $n$ and $m$ lies in detecting enough algebraic cycles to cover all Hodge classes. For this reason it is natural to look for some other families of varieties, symmetric enough to have a refined spectral decomposition on their cohomologies, but with a simpler subgroup of automorphisms, so one can expect to have better control over the algebraic cycles we can produce on them. In this article we study the Hodge cycles of Klein varieties, which correspond to
$$
Y^n_m:=\{K:=y_0^{m-1}y_1+y_1^{m-1}y_2+\cdots+y_n^{m-1}y_{n+1}+y_{n+1}^{m-1}y_0=0\}\subseteq\P^{n+1}.
$$
Such a variety admits the following diagonal automorphism of order $d:=\frac{(m-1)^{n+2}-(-1)^{n+2}}{m}$ (see for instance \cite{GLMV22})
$$
\sigma(y_0:y_1:\cdots:y_{n+1}):=(\zeta_{md}y_0:\zeta_{md}^{1-m}y_1:\cdots:\zeta_{md}^{(1-m)^{n+1}}y_{n+1}).
$$
Since $\sigma$ has order $d$, its induced action on $\mathrm{H}^n(Y^n_m,\Q)_\prim$ has order a divisor of $d$, hence all its eigenvalues are $d$-roots of unity. We denote the eigenspace of an eigenvalue $\zeta_d^i$ by
$$
V(i):=\{\omega\in \mathrm{H}^n(Y^n_m,\C)_\prim: \sigma^*(\omega)=\zeta_d^i\cdot\omega\}
$$
for every $i\in\Z/d$. We have the following result analogous to Shioda's decomposition.

\begin{thm}\label{thm0.1}
The induced action of $\sigma$ on $\mathrm{H}^n(Y^n_m,\C)_\prim$ is faithful and induces the following spectral decomposition
\begin{equation}
\label{eq1}
\mathrm{H}^n(Y^n_m,\C)_\prim=\bigoplus_{i=0}^{d-1} V(i)
\end{equation}
satisfying the following properties
\begin{itemize}
    \item[(i)] $\dim_\C V(i)=1$ for every $i=1,\ldots,d-1$,
    \item[(ii)] $\dim_\C V(0)=1+(-1)^n$ and it is generated by linear cycles,
    \item[(iii)] for $p+q=n$ and $i\in\{1,\ldots,d-1\}$ let $\alpha(i)=(a_0,\ldots,a_{n+1})\in\{1,\ldots,m-1\}^{n+2}$ such that $mi\equiv_{md}\sum_{j}a_j(1-m)^j$, then
    $$
\mathrm{H}^{p,q}(Y^n_m)_\prim=\left\{\begin{array}{cc}
   \bigoplus_{|\alpha(i)|=q+1}V(i)  & \text{ if } 2p\neq n,\\
   \bigoplus_{|\alpha(i)|=p+1}V(i)\oplus V(0) & \text{ if }2p=n.
\end{array}\right.
    $$
    \item[(iv)] for $n=2p$ even
    $$
\mathrm{Hod}^n(Y^n_m)_\prim\otimes\C=\bigoplus_{i\in\mathfrak{B}^n_m(\cK)}V(i)\oplus V(0),
    $$
    where $\fB^n_m(\cK):=\{i\in\{1,\ldots,d-1\}: |\alpha(t\cdot i)|=p+1\ , \ \forall t\in(\Z/d)^\times\}$.
\end{itemize}
\end{thm}
Klein varieties are also examples of Delsarte hypersurfaces, which are smooth hypersurfaces of $\P^{n+1}$ defined by a sum of $n+2$ monomials whose exponent matrix is invertible. In \cite{shioda1986explicit}, Shioda gave a general algorithm to study the N\'eron-Severi group of Delsarte surfaces, by means of a rational monomial covering from a Fermat surface, under which the transcendental part of the Hodge structure of the Delsarte surface becomes a sub-Hodge structure of the Fermat surface. These monomial coverings also exist in higher dimensions. For Klein hypersurfaces the monomial covering is constructed as follows. Let $A$ be the exponent matrix of $K$ and $B=\text{adj}(A)$, note that $\det (A)=(m-1)^{n+2}-1=md$. Then the monomial covering
$\phi\colon X_{md}^n\dashrightarrow Y_m^n$, from the degree $md$ Fermat variety to the Klein variety is given in coordinates by $$y_i=\prod_j x_j^{B_{ij}}.$$
We use this monomial covering to relate the Hodge cycles of the Klein variety to Hodge cycles of the Fermat variety. Solving the indeterminacy of $\phi$ we get the diagram
\[
\begin{tikzcd}[column sep=small]
 & \widetilde{X} \arrow[dl, "\varepsilon"'] \arrow[dr, "\widetilde{\phi}"] & \\
X^n_{md} \arrow[rr, dashrightarrow, "\phi"] & & Y^n_m
\end{tikzcd}
\]
and we show (\cref{prop2.1}) that for each $i=1,\ldots,d-1$ 
$$
\tilde\phi^*(V(i))=\varepsilon^*(V(\alpha(i)\cdot B)).
$$
This reduces the Hodge conjecture on the Klein variety $Y^n_m$ to studying the set of characters 
$$
\mathfrak{K}^n_{md}:=\{\alpha(i)\cdot B\in\mathfrak{B}^n_{md}:i\in\mathfrak{B}^n_m(\cK)\}
$$
and determining when $\mathfrak{K}^n_{md}\subseteq\mathfrak{C}^n_{md}$. In particular, it holds whenever $X^n_{md}$ satisfies the Hodge conjecture. In analogy with Shioda's notation for characters in Fermat varieties, we separate the set $\fB^n_m(\cK)$ of Hodge characters of $Y^n_m$ into two parts. The first part corresponds to characters supported in rational subvarieties obtained as images of linear subvarieties under $\phi$
$$
\fD^n_m(\cK):=\{i\in\fB^n_m(\cK): \alpha(i)\cdot B\in\fD^n_{md}\}.
$$
The remaining part is Galois invariant and so it can be partitioned in orbits
$$
\fB^n_m(\cK)\setminus\fD^n_m(\cK)=\bigsqcup_{e\in \mathcal{E}^n_m(\cK)}\frac{d}{e}\cdot(\Z/e)^\times
$$
where $\mathcal{E}^n_m(\cK):=\{e\mid d:\frac{d}{e}\in \fB^n_m(\cK)\setminus\fD^n_m(\cK)\}$ are called \textit{exceptional denominators} and the corresponding characters $\frac{d}{e}\cdot(\Z/e)^\times\subseteq\fB^n_m(\cK)\setminus\fD^n_m(\cK)$ are called \textit{exceptional Hodge characters}.
\begin{thm}\label{thm0.2}
Let $\mathcal{E}^n(\cK):=\bigcup_{m\in\N}\mathcal{E}^n_m(\cK)$. Then $\mathcal{E}^n(\cK)$ is a finite set. Furthermore, for every $e\in\mathcal{E}^n(\cK)$ there is a set $S_e\subseteq\Z/e$ such that for every $m\in\N$
$$
\mathcal{E}^n_m(\cK)=\{e\in\mathcal{E}^n(\cK): [m]_e\in S_e\}.
$$
In particular, $\mathcal{E}^n_m(\cK)$ is completely determined by the residue class of $m$ modulo $\gcd(\mathcal{E}^n(\cK)).$
\end{thm}

For Fermat surfaces, Shioda \cite{Shioda1981} also introduced a set of exceptional characters, which was later proved by Aoki \cite[Theorem C]{Aoki1983} to be a finite set (consisting of characters appearing in degree at most $180$). The finiteness of this set is what guarantees the termination of Shioda's algorithm \cite{shioda1986explicit} for computing the Picard rank of Delsarte surfaces. For higher dimension Fermat varieties, it is not clear how to generalize Shioda's notion of exceptional characters in such a way that it remains being finite for a fixed dimension (see \cite[Conjecture 4.1]{MMRV2026} for a possible direction). In this sense,  \cref{thm0.2} has the advantage of showing the finiteness of exceptional characters for all dimensions. The proof of \cref{thm0.2} is based on a fine analysis of odd Dirichlet characters of maximal conductor in an admissible set (see \cref{thmAIV} and its proof). The relation of Hodge characters and Dirichlet character follows from a result of Aoki \cite[Proposition 2.2]{Aoki1983}. As a consequence of that analysis we can provide sufficient conditions to the Hodge conjecture for Klein varieties.

\begin{thm}\label{thm0.3}
Let $m,n\in\N$ with $n$ even. Assume that for all primes $p$ such that $p-1\mid n+2$ it holds that $p\nmid d$. Then $\fB^n_m(\cK)=\fD^n_m(\cK)$.
\end{thm}

We remark that by \cite[Theorem A]{Aoki1983} it follows that $\fB_m^n(\cK)=\fD^{n}_m(\cK)$ when $\gcd(md,(n+1)!)=1$. Due to the fact $\fB^n_{md}=\fD^n_{md}$ in that case. \cref{thm0.3} covers the case $\gcd(md,(n+3)!)=1$, but it also covers infinitely many new cases where $\fB_m^n(\cK)=\fD^{n}_m(\cK)$ but $\fB^n_{md}\neq\fD^n_{md}$ (essentially what happens is $\mathfrak{K}^n_{md}\subseteq\fD^n_{md}\subsetneq \fB^n_{md}$). 

Once we know the finiteness of exceptional characters, it is natural to try to determine them. However, the upper bound provided by the proof is too big to perform a complete computer search in a fixed dimension. For small dimensions $n=2, 4, 6, 8$ we perform a computer based search of the exceptional denominators up to degree $m\le 10^6$ and obtain the following result. 

\begin{thm}\label{thm0.4}
For $n=2,4,6$ and $m\leq 10^7$, and for $n=8$ and $m\le 10^6$ we have the following equalities
\[
\cE_m^2(\cK)=\Big\{20_{\,([m]_{20}=4,8)},\ 30_{\,([m]_{30}=8,14)},\
60_{\,([m]_{60}=8,14,38,44)},\ 120_{\,([m]_{120}=14,38)}\Big\},
\]
\[
\renewcommand{\arraystretch}{1.5}
\begin{array}{rcl}
\cE_m^4(\cK)&=&\Big\{14_{([m]_{14}=4,6)},\ 21_{([m]_{21}=3,12)},\
 28_{([m]_{28}=6,18)},\ 39_{([m]_{39}=5,11)},\\
&&\ \ 42_{([m]_{42}=6,12,18,20,24,32)},\
 78_{([m]_{78}=18,24,30,36,44,50)},\ 84_{([m]_{84}=6,18)},\\
&&\ \ 126_{([m]_{126}=32,62)},\ 156_{([m]_{156}=18,102)},\
 168_{([m]_{168}=18,90)},\ 210_{([m]_{210}=102,132)}\Big\},\\
\cE_m^6(\cK)&=&\Big\{20_{([m]_{20}=4,8)},\ 30_{([m]_{30}=8,14)},\
 60_{([m]_{60}=8,14,38,44)},\ 102_{([m]_{102}=20,26,44,50)},\\
&&\ \ 120_{([m]_{120}=14,38)},\ 170_{([m]_{170}=54,78,84,128)},\
 204_{([m]_{204}=20,44,128,152)},\\
&&\ \ 510_{([m]_{510}=128,224,248,254)}\Big\},\\
\cE_m^8(\cK)&=&\Big\{33_{([m]_{33}=8,14,20,29)},\ 55_{([m]_{55}=5,10,15,50)},\
 62_{([m]_{62}=16,24,28,30)},\\
&&\ \ 66_{([m]_{66}=8,14,20,62)},\
 110_{([m]_{110}=10,20,30,40,42,50,52,60,62,70,80,102)},\\
&&\ \ 124_{([m]_{124}=30,78,86,90)},\ 186_{([m]_{186}=86,92,140,152)},\
 220_{([m]_{220}=30,130,150,190)},\\
&&\ \ 330_{([m]_{330}=20,80,140,260)},\ 363_{([m]_{363}=41,95,113,119)},\
 682_{([m]_{682}=214,334,340,512)},\\
&&\ \ 726_{([m]_{726}=404,458,476,482)},\
 2046_{([m]_{2046}=512,896,1016,1022)}\Big\}.
\end{array}
\]
\end{thm}

Since for Fermat surfaces it is known that all exceptional characters appear up to degree $180$, the first equality holds for any $m\in\N$. It is also possible to perform similar computations for dimension $10$, but due to their size we decided not to include it in the previous result (cf. \cite{tesisMario}). As an immediate consequence of \cref{thm0.4} and \cite[Theorem 1.1]{MMRV2026} we obtain more cases of the Hodge conjecture for Klein varieties of small dimension, stated as follows.

\begin{cor}\label{cor0.5}
The Hodge conjecture for $Y^n_m$ holds for $n=2,4$. For $n=6$ and $m\le 10^7$ it holds if $[m]_{204}\neq20,44,128, 152$. For $n=8$ and $m\le 10^6$ it holds if $[m]_{220}\neq30,130,150,190$ and $[m]_{2046}\neq512, 896,1016,1022$.
\end{cor}

The article is organized as follows. In \cref{sec1} we prove \cref{thm0.1} (see \cref{prop1.1} and \cref{prop1.2}). In \cref{sec2} we relate the characters of the Klein varieties to characters of the covering Fermat variety (\cref{prop2.1}). Using this correspondence we give an explicit description of $\fD^n_m(\cK)$ (\cref{prop2.2}). In \cref{sec3} we study general properties of the set of exceptional denominators (\cref{thmAIV}) and prove \cref{thm0.2} and \cref{thm0.3}. Finally in \cref{sec4} we characterize the set of exceptional denominators for small dimensions and degrees, providing the proof of \cref{thm0.4} and \cref{cor0.5}.

\bigskip 

\noindent\textbf{Acknowledgments and AI disclosure.} The authors are grateful to Maximiliano Miranda for his support on some computational verifications of the article. The first author was supported by the program PIIC N$^\circ$ 015/2025 of the Postgraduate Direction of UTFSM. The second author was supported by Fondecyt ANID regular grant 1240101 and Fondecyt ANID initiation grant 11251404. It was used the help of Claude Opus 5.5 to complete the proof of \cref{thm0.2}, which was then studied, verified, simplified and written by the authors. All the remaining mathematical results and the writing of the article are integrally due to the authors. The authors take full responsibility for all the mathematical statements and proofs in the article.

\section{Hodge characters on Klein varieties}
\label{sec1}
In this section we study the spectral decomposition for Klein varieties, providing a proof of \cref{thm0.1}. Since the Klein variety $Y^n_m$ is a smooth hypersurface of $\P^{n+1}$ given by $Y^n_m=\{K=0\}$, it is well known that Griffiths' basis theorem \cite[Corollary 6.12]{vo03} identifies the graded pieces of the Hodge filtration of primitive cohomology of $Y^n_m$ with certain homogeneous parts of the Jacobian ring of $K$, via residue forms as follows: For $p+q=n$ and $\Omega=\iota_E(dy_0\wedge\cdots\wedge dy_{n+1})=\sum_{i=0}^{n+1}(-1)^iy_i\widehat{dy_i}$ we have the isomorphism
$$
P\in J(K)_{m(q+1)-n-2}\xrightarrow{\sim} \text{res}\left(\frac{P\Omega}{K^{q+1}}\right)^{p,q}\in \mathrm{H}^{p,q}(Y^n_m)_\prim.
$$
\begin{prop}\label{prop1.1}
For $p+q=n$ and $i\in\{1,\ldots,d-1\}$ let $\alpha(i)=(a_0,\ldots,a_{n+1})\in\{1,\ldots,m-1\}^{n+2}$ such that $mi\equiv_{md}\sum_{j}a_j(1-m)^j$, then
    $$
\mathrm{H}^{p,q}(Y^n_m)_\prim=\left\{\begin{array}{cc}
   \bigoplus_{|\alpha(i)|=q+1}V(i)  & \text{ if } 2p\neq n,\\
   \bigoplus_{|\alpha(i)|=p+1}V(i)\oplus V(0) & \text{ if }2p=n.
\end{array}\right.
    $$
Furthermore, for $i\in\{1,\ldots,d-1\}$, $$V(i)=\C\cdot\omega_{\alpha(i)}$$ is generated by
$$
\omega_{\alpha(i)}:=\text{res}\left(\frac{y^{\alpha(i)-\bone}\Omega}{K^{q+1}}\right)\in \mathrm{H}^{p,q}(Y^n_m)_\prim.
$$
And when $n=2p$, $$V(0)=\C\cdot\omega_0\oplus\C\cdot\omega_1$$ is generated by the classes $$\omega_0:=\text{res}\left(\frac{(y_0y_2\cdots y_{n})^{m-2}\Omega}{K^{p+1}}\right)\in \mathrm{H}^{p,p}(Y^n_m)_\prim,
$$
$$
\omega_1:=\text{res}\left(\frac{(y_1y_3\cdots y_{n+1})^{m-2}\Omega}{K^{p+1}}\right)\in \mathrm{H}^{p,p}(Y^n_m)_\prim.$$
\end{prop}

\begin{proof}
From Grifiths' basis we know $\mathrm{H}^{p,q}(Y^n_m)_\prim$ is generated by a basis of terms of the form $\omega_\alpha:=\text{res}\left(\frac{y^{\alpha-\bone}\Omega}{K^{q+1}}\right)^{p,q}$ with $\alpha-\bone\in\{0,\ldots,m-2\}^{n+2}$ and $|\alpha|=q+1$. For such a class $\alpha=(a_0,\ldots,a_{n+1})\in\{1,\ldots,m-1\}^{n+2}$ we have
$$
\sigma^*\omega_\alpha=\zeta_{md}^{(\sum_{j=0}^{n+1}(a_j-1)(1-m)^j)-d}\omega_\alpha=\zeta_{md}^{\sum_{j=0}^{n+1}a_j(1-m)^j}\omega_\alpha=\zeta_d^i\omega_\alpha
$$
for $mi\equiv_{md}\sum_ja_j(1-m)^j$. Note that since $|\alpha|=q+1$ it follows that $\sum_ja_j(1-m)^j\equiv_m0$ and so the map $\alpha\mapsto i\in\{0,1,\ldots,d-1\}$ is well defined. On the other hand, for $i\in\{0,\ldots,d\}$ we can define $\alpha(i)$ which is given by writing $mi=\sum_{j=0}^{n+1}b_j(m-1)^j$ in its $(m-1)$-adic expansion and then take $a_j:=b_j+1$ for $j$ even and $a_j:=m-1-b_j$ for $j$ odd. It is not hard to see that this last map is a bijection and that $\omega_{\alpha(i)}\in V(i)$ for all $i=1,\ldots,d-1$, and that $\omega_{\alpha(0)},\omega_{\alpha(d)}\in V(0)$. From this, all statements follow. In particular, for each $\alpha$, since $\text{res}\left(\frac{y^{\alpha-\bone}\Omega}{K^{q+1}}\right)\in V(i)\subseteq \mathrm{H}^{p,q}(Y^n_m)_\prim$, all these residue forms are pure of type $(p,q)$. Note that $\omega_{\alpha(0)}=\omega_1$ and $\omega_{\alpha(d)}=\omega_0$. 
\end{proof}

Since the spectral decomposition has eigenvalues in the cyclotomic field $\Q(\zeta_d)$, the decomposition also holds over this field. In fact, the eigenspace $V(i)$ is defined over $\Q(\zeta_d^i)$. In particular, $V(0)$ is defined over $\Q$. For any subfield $K\subseteq\C$, let us denote by
$$
V(i)_{K}:=\mathrm{H}^n(Y^n_m,K)_\prim\cap V(i),
$$
then 
$$
\mathrm{H}^n(Y^n_m,\Q(\zeta_d))_\prim=\bigoplus_{i=0}^{d-1}V(i)_{\Q(\zeta_d)}.
$$
Furthermore, the Galois group $\Gamma=\text{Gal}(\Q(\zeta_d)/\Q)\simeq(\Z/d)^\times$ acts on this space in terms of periods (see \cite[\S 2.5]{MMRV2026}) as follows: for $\tau\in\Gamma$, $\omega\in \mathrm{H}^n(Y^n_m,\Q(\zeta_d))$ and $\delta\in \mathrm{H}_n(Y^n_m,\Z)$
$$
\int_\delta\tau(\omega)=\tau\left(\int_\delta\omega\right).
$$
In particular, the space of primitive classes is Galois invariant.

\begin{prop}\label{prop1.2}
Let $\tau\in\Gamma$ be such that $\tau(\zeta_d)=\zeta_d^t$ for $t\in(\Z/d)^\times$, and $i\in\{0,\ldots,d-1\}$. Then
$$
\tau(V(i)_{\Q(\zeta_d)})=V(t\cdot i)_{\Q(\zeta_d)}.
$$
In particular, for $n=2p$ even
$$
\mathrm{Hod}^n(Y^n_m)_\prim\otimes\C=\bigoplus_{i\in\fB^n_m(\cK)}V(i)\oplus V(0),
$$
where $\fB^n_m(\cK):=\{i\in\{1,\ldots,d-1\}: |\alpha(t\cdot i)|=p+1\ , \ \forall t\in(\Z/d)^\times\}$. Moreover,
$$
V(0)_\Q=\Q\cdot[L_0]_\prim\oplus\Q\cdot[L_1]_\prim,
$$
where $L_0:=\{y_0=y_2=\cdots=y_n=0\}$ and $L_1=\{y_1=y_3=\cdots=y_{n+1}=0\}$.
\end{prop}

\begin{proof}
For $i\in\{0,1,\ldots,d\}$, recall from the proof of \cref{prop1.1} that $\alpha(i)=(a_0,\ldots,a_{n+1})\in\{1,\ldots,m-1\}^{n+2}$ is constructed such that
$$
\sigma^*\omega_{\alpha(i)}=\zeta_d^i\omega_{\alpha(i)}.
$$
For $i\equiv_d 0$, it is clear that $V(0)_{\Q(\zeta_d)}$ is Galois invariant, since it is defined over $\Q$. For $i\in\{1,\ldots,d-1\}$, we know by \cref{prop1.1} that $V(i)=\C\cdot\omega_{\alpha(i)}$, and so 
$$
V(i)_{\Q(\zeta_d)}=\Q(\zeta_d)\cdot\eta_{\alpha(i)}
$$
for $\eta_{\alpha(i)}=c_{\alpha(i)}\cdot \omega_{\alpha(i)}$ for some $c_{\alpha(i)}\in\C^\times$. In consequence, $\sigma^*\eta_{\alpha(i)}=\zeta_d^i\eta_{\alpha(i)}$.
Then, for all $\delta\in \mathrm{H}_n(Y^n_m,\Z)$
$$
\int_\delta\sigma^*(\tau(\eta_{\alpha(i)}))=\int_{\sigma_*(\delta)}\tau(\eta_{\alpha(i)})=\tau\left(\int_{\sigma_*(\delta)}\eta_{\alpha(i)}\right)=\tau\left(\int_\delta\sigma^*\eta_{\alpha(i)}\right)=\int_\delta\tau(\zeta_d^i)\cdot\tau(\eta_{\alpha(i)}),
$$
and so $\tau(\eta_{\alpha(i)})\in V(t\cdot i)_{\Q(\zeta_d)}$ for $\tau(\zeta_d)=\zeta_d^t$. Since these spaces are one dimensional the first part of the proposition follows. The second part of the proposition follows directly from this. So, we are left with the last assertion. Assume $n=2p$ is even, then by \cite[Theorem 1.1]{Vil22PCIC} it is a direct computation to get that
$$
\C\cdot[L_0]_\prim=\C\cdot(\omega_0+(-1)^p(m-1)^{p+1}\omega_1),
$$
$$
\C\cdot[L_1]_\prim=\C\cdot((-1)^p(m-1)^{p+1}\omega_0+\omega_1).
$$
In particular, $V(0)$ is generated by $[L_0]_\prim$ and $[L_1]_\prim$, and since these are rational classes, they generate $V(0)_\Q$.
\end{proof}

\begin{rmk}
The previous proposition, together with \cref{prop2.1} complete the proof of \cref{thm0.1}.
\end{rmk}

\section{Monomial covering}
\label{sec2}
Let us now relate the Hodge structure and the spectral decomposition of $Y^n_m=\{K=y_0^{m-1}y_1+\cdots+y_{n+1}^{m-1}y_0=0\}$ with the one of a Fermat variety lying over it under a rational monomial covering $$\phi:X^n_{md}\dashrightarrow Y_m^n.$$
The construction of this monomial map is classical \cite{Shioda1981,shioda1983algebraic} and goes as follows: let
$A$ be the exponent matrix of $K$, $\det(A)=(m-1)^{n+2}-1=md$ and 
$$
B=\text{adj}(A)=\begin{pmatrix}(m-1)^{n+1}&-(m-1)^n&(m-1)^{n-1}&\cdots&-1\\-1&(m-1)^{n+1}&-(m-1)^n&\cdots&(m-1)\\\vdots&\vdots&\vdots&\ddots&\vdots\\-(m-1)^n&(m-1)^{n-1}&-(m-1)^{n-2}&\cdots&(m-1)^{n+1}\end{pmatrix}.
$$
For $X^n_{md}=\{x_0^{md}+\cdots+x_{n+1}^{md}=0\}$, the monomial map is given in coordinates by
$$\phi(x)_i=\prod_j x_j^{B_{ij}}.$$
Solving the indeterminacy of $\phi$ we get the diagram
\[
\begin{tikzcd}[column sep=small]
 & \widetilde{X} \arrow[dl, "\varepsilon"'] \arrow[dr, "\widetilde{\phi}"] & \\
X^n_{md} \arrow[rr, dashrightarrow, "\phi"] & & Y^n_m .
\end{tikzcd}
\]
\begin{prop}\label{prop2.1}
For each $i=1,\ldots,d-1$ 
$$
\tilde\phi^*(V(i))=\varepsilon^*(V(\alpha(i)\cdot B)).
$$
\end{prop}

\begin{proof}
Since $V(i)=\C\cdot\omega_{\alpha(i)}$, it is enough for us to compute $\tilde{\phi}^*\omega_{\alpha(i)}$. Let $\beta=\alpha(i)-\bone$. Clearly $\phi(x)^\beta=x^{\beta B}$, so it remains to compute $\phi^*\Omega(y)$.
Write $B_k=(B_{k,0},\dots,B_{k,n+1})$ for the $k$-th row of $B$, so that
$\phi(x)_k=x^{B_k}$ and
\[
d\big(\phi(x)_k\big)=x^{B_k}\sum_{j=0}^{n+1}\frac{B_{k,j}}{x_j}\,dx_j .
\]
Substituting into
$\Omega(y)=\sum_k(-1)^ky_k\,dy_0\wedge\cdots\widehat{dy_k}\cdots\wedge dy_{n+1}$
and factoring out the monomials gives
\[
\Omega(\phi(x))=\sum_{k=0}^{n+1}(-1)^k\Big(\prod_{r}x^{B_r}\Big)
\bigwedge_{r\neq k}\Big(\sum_{j}B_{r,j}\frac{dx_j}{x_j}\Big).
\]
Expanding the wedge by multilinearity,
\[
\bigwedge_{r\neq k}\Big(\sum_{j}B_{r,j}\frac{dx_j}{x_j}\Big)
=\sum_{j=0}^{n+1}\det\big(B^{(k,j)}\big)\,
\frac{dx_0}{x_0}\wedge\cdots\widehat{\frac{dx_j}{x_j}}\cdots
\wedge\frac{dx_{n+1}}{x_{n+1}} ,
\]
where $B^{(k,j)}$ is the minor obtained by deleting row $k$ and column $j$. Since
every column of $B$ sums to $d$ we have
$\prod_r x^{B_r}=x^{\,d\bone}$, so
\[
\Omega(\phi(x))=x^{d\bone}\sum_{j=0}^{n+1}
\Big(\sum_{k=0}^{n+1}(-1)^k\det\big(B^{(k,j)}\big)\Big)
\bigwedge_{r\neq j}\frac{dx_r}{x_r} .
\]
Cofactor expansion together with the constancy of the column sums yields
\[
\sum_{k=0}^{n+1}(-1)^k\det\big(B^{(k,j)}\big)=(-1)^j\,\frac{\det B}{d}
=(-1)^j\,\frac{(md)^{n+2}}{d} ,
\]
and therefore
\[
\Omega(\phi(x))=\frac{(md)^{n+2}}{d}\,x^{\,d\bone-\bone}
\sum_{j}(-1)^jx_j\bigwedge_{r\neq j}dx_r
=\frac{(md)^{n+2}}{d}\,x^{\,d\bone-\bone}\,\Omega(x).
\]
Since $K(\phi(x))=F(x)$, applying the residue map gives
\[
\phi^{*}(\omega_{\alpha(i)})
=\phi^{*}\res\!\left(\frac{y^\beta\Omega(y)}{K(y)^{q+1}}\right)^{\!p,q}
=\frac{(md)^{n+2}}{d}\,\res\!\left(
\frac{x^{\,\beta B+d\bone-\bone}\,\Omega(x)}{F(x)^{q+1}}\right)^{\!p,q}
=\frac{(md)^{n+2}}{d}\,\omega_{\alpha(i)\cdot B} .
\]
Taking $\varepsilon^*$ we obtain the result. Note that for $i\in\{0,\ldots,d\}$, $\alpha(i)\cdot B$ has non-zero entries if and only if $\alpha(i)\neq(1,m-1,\ldots,1,m-1), (m-1,1,\ldots,m-1,1)$, which happens if and only if $i\neq 0,d$.
\end{proof}

In consequence, the Hodge conjecture for $Y^n_m$ is reduced to analyze the algebraicity of the characters in
$$
\mathfrak{K}^n_{md}:=\{\alpha(i)\cdot B\in\mathfrak{B}^n_{md}:i\in\mathfrak{B}^n_m(\cK)\}.
$$
In fact, $V(0)$ was already shown to be generated by algebraic cycles, and if $Z\subseteq X^n_{dm}$ is an algebraic subvariety such that $[Z]_\prim$ has a non-trivial projection on $V(\alpha(i)\cdot B)$, then $[\tilde\phi_*\varepsilon^*(Z)]_\prim$ has a non-trivial projection on $V(i)$.

\begin{prop}\label{prop2.2}
For each $i\in\{1,\ldots,d-1\}$, 
$$
\alpha(i)\cdot B=(-mi(1-m)^{n+1},-mi(1-m)^n,\ldots,-mi(1-m),-mi)\in\mathfrak{A}^n_{md}.
$$
In particular, $i\in\fB^n_m(\cK)$ if and only if $\alpha(i)\cdot B\in \fB^n_{md}$ if and only if 
$$
\alpha_0:=(1,1-m,(1-m)^2,\ldots,(1-m)^{n+1})\in\fB^n_e
$$
for $e=d/\gcd(i,d)$. And we have similar equivalences for algebraic characters. Furthermore, $i\in\fD^n_m(\cK)$ if and only if there exists some $h\mid \frac{n}{2}+1$, such that
$$
e\mid{(1-m)^h+1} .
$$
\end{prop}

\begin{proof}
For the last column of $B$ we have
$B_{i,n+1}=-(1-m)^{i}$ for $i=0,\ldots,n+1$, hence
\[
\big(\alpha(i)\cdot B\big)_{n+1}
=-\!\!\sum_{k}\!\!a_k(1-m)^k\equiv_{md}-mi.
\]
Since each column of $B$ is
$(1-m)$ times the next modulo $md$, the coordinates of $\alpha(i)\cdot B$ satisfy
$\big(\alpha(i)\cdot B\big)_{j-1}=(1-m)\big(\alpha(i)\cdot B\big)_{j}$, and reading
backwards from $j=n+1$ gives
\[
\alpha(i)\cdot B\equiv_{md} (-mi(1-m)^{n+1},-mi(1-m)^n,\ldots,-mi(1-m),-mi).
\]
Hence, using the Galois action of $\text{Gal}(\Q(\zeta_{md})/\Q)\simeq(\Z/md)^\times$ on the characters of the Fermat variety, we have that $\alpha(i)\cdot B\in\fB^n_{md}$ if and only if $((1-m)^{n+1}md/e,(1-m)^nmd/e,\ldots,(1-m)md/e,md/e)\in\fB^n_{md}$ for $d/e=\gcd(i,d)$. By \cite[Proposition 3.2]{MMRV2026} this holds if and only if $((1-m)^{n+1},\ldots,1-m,1)\in\fB^n_{e}$, and by \cite[Proposition 3.4]{MMRV2026} this happens if and only if $\alpha_0\in\fB^n_{e}$. The same propositions, together with \cite[Remark 3.2]{MMRV2026}, imply the corresponding statement for algebraic characters. For the last statement, it is clear that if $i\in\fD^n_m(\cK)$, then $\alpha_0\in\fD^n_{e}$ and so ${e}\mid 1+(1-m)^k$ for some $0<k<n+2$. Assume $h$ is the smallest among those $k$. If $
{e}=2$, then $h=1$ and so $2h=2\mid n+2$. If ${e}>2$, then $\text{ord}_{{e}}(1-m)=2h$. Since $(1-m)^{n+2}\equiv_{{e}}1$, it follows that $2h\mid n+2$. Conversely, if the condition holds, then clearly $\alpha_0=(1,1-m,\ldots,(1-m)^{h-1},-1,-(1-m),\ldots,-(1-m)^{h-1},\ldots)\in\fD^n_{e}$ and so $\alpha(i)\cdot B\in\fD^n_{md}$.
\end{proof}

\begin{rmk}\label{rmk2.1}
Because of the previous proposition, to determine whether $i\in\fB^n_m(\cK)$ (or $\fC^n_m(\cK)$, $\fD^n_m(\cK)$) is equivalent to determine whether $\alpha_0\in\fB^n_{e}$ (respectively $\fC^n_{e}$, $\fD^n_{e}$). For this reason we say $\alpha_0$ is the test vector.
\end{rmk}

\section{Exceptional denominators}
\label{sec3}
In the previous section we described the Hodge characters of $Y^n_m$ supported in rational subvarieties $\fD^n_m(\cK)$. Since the complement is Galois invariant, we partition it into orbits as
\begin{equation}\label{eq3.1}
\fB^n_m(\cK)\setminus\fD^n_m(\cK)=\bigsqcup_{e\in \mathcal{E}^n_m(\cK)}\frac{d}{e}\cdot(\Z/e)^\times    
\end{equation}
where $\mathcal{E}^n_m(\cK):=\{e\mid d: \frac{d}{e}\in \fB^n_m(\cK)\setminus\fD^n_m(\cK)\}$. We say the elements in $\mathcal{E}^n_m(\cK)$ are the \textit{exceptional denominators} for $Y^n_m$. In this section we study the structure of these numbers.

\begin{rmk}\label{rmk3.1}
By \cref{prop2.2} we see that whether some number $e\in\N$ belongs to $\mathcal{E}^n_m(\cK)$ is equivalent to have that $e\mid d$, $\alpha_0\in\fB^n_e$ and $e\nmid (1-m)^h+1$ for all $h\mid \frac{n}{2}+1$. Since these are all properties of $m$ in modulo $e$, we see that this is determined by the residue class $[m]_e\in\Z/e$. This leads us to define the following.
\end{rmk}

\begin{dfn}
Let $n\in\N$. We define $\cE^n(\cK):=\bigcup_{m\in\N}\cE^n_m(\cK)$. For every $e\in\cE^n(\cK)$ let 
$$
S_e:=\{[m]_e\in\Z/e: e\in\cE^n_m(\cK)\}=\{[m]_e\in\Z/e: \alpha_0=(1,1-m,\ldots,(1-m)^{n+1})\in\fB^n_e\}.
$$
\end{dfn}

As a consequence of the previous remark we have that
$$
\cE^n_m(\cK)=\{e\in\cE^n(\cK):[m]_e\in S_e\}.
$$
And so we reduce the study of exceptional denominators to understand the sets $\cE^n(\cK)$.  

\begin{rmk}
The key idea to prove the finiteness of $\cE^n(\cK)$ is to use the a characterization, due to Aoki, of Hodge characters of Fermat varieties in terms of Dirichlet characters. Let us first recall some standard notations about Dirichlet characters. For $m\in\Z$ we denote by $C^-(m)$ the set of odd characters modulo $m$, i.e. characters $\chi:(\Z/m)^\times\to\C^\times$ such that $\chi(-1)=-1$. We denote by $PC(m)$ the set of primitive characters modulo $m$, i.e. those which can not be factored by $(\Z/m)^\times\to(\Z/m')^\times$ for some divisor $m'\mid m$. When $\chi\in PC(m)$ we say $\chi$ has conductor $\text{cond}(\chi)=m$. In general, for a character $\chi:(\Z/m)^\times\to\C^\times$, we say it has conductor $\text{cond}(\chi)=m'$ if $m'\mid m$ and $\chi=\chi'\circ\pi$ for $\pi:(\Z/m)^\times\to(\Z/m')^\times$ and $\chi'\in PC(m')$. By the Chinese remainder theorem we can always write in a unique way $$\chi=\prod_{p\text{ prime}}\chi_p$$ where $\chi_p$ is a character modulo $p^{v_p}\mid\mid\text{cond}(\chi)$. In particular, $\chi_p=1$ is trivial if and only if $p\nmid \text{cond}(\chi)$. Furthermore, $\chi$ is primitive if and only if all non-trivial factors $\chi_p$ are primitive. We denote by $PC^-(m)=PC(m)\cap C^-(m)$ the primitive odd characters modulo $m$.
\end{rmk}

\begin{prop}\label{prop3.1}
Let $e\mid d=\frac{(m-1)^{n+2}-1}{m}$ and $\alpha_0=(1,1-m,\ldots,(1-m)^{n+1})\in((\Z/e)^\times)^{n+2}$. Then $\alpha_0\in\fB^n_e$ if and only if for every $f\mid e$ and every $\chi\in PC^-(e/f)$ we have
$$
\chi(1-m)\neq 1 \ \ \text{ or } \ \ \exists \text{ prime } p\mid f,\  p\nmid e/f: \chi(p)=1. 
$$
\end{prop}

\begin{proof}
Note that since $\chi(1-m)^{n+2}=\chi((1-m)^{n+2})=\chi(1)=1$, then
$$
\chi(\alpha_0)=\sum_{i=0}^{n+1}\chi((1-m)^i)=\sum_{i=0}^{n+1}\chi(1-m)^i=\left\{\begin{array}{cc}
    n+2 & \text{if }\chi(1-m)=1, \\
    0 & \text{if }\chi(1-m)\neq 1. 
\end{array}\right.
$$
Hence, the result follows by \cite[Proposition 2.2]{Aoki1983}.
\end{proof}

Before the main theorem of the section let us prove a short elementary lemma.

\begin{lemma}\label{lemma}
Let $m,n\in\N$ both even, and put $a:=v_2(n+2)$. Then, for $d=\frac{(m-1)^{n+2}-1}{m}$ we have
\[
  v_2(d)=v_2(m-2)+a-1 .
\]
\end{lemma}

\begin{proof}
Write $n+2=2^{a}j$ with $j$ odd and factor
\[
  (m-1)^{n+2}-1=((m-1)^{j}-1)\,((m-1)^{j}+1)\prod_{i=1}^{a-1}\bigl((m-1)^{2^{i}j}+1\bigr).
\]
Since $j$ is odd, $(m-1)^{j}-1=(m-2)((m-1)^{j-1}+\cdots+(m-1)+1)$ and the second factor is a
sum of $j$ odd terms, hence odd; therefore $v_2((m-1)^{j}-1)=v_2(m-2)$, and the same
argument gives $v_2((m-1)^{j}+1)=v_2(m)$. For $i\ge 1$, $(m-1)^{2^{i}j}$ is an odd
square, so $(m-1)^{2^{i}j}\equiv 1\pmod 8$ and $(m-1)^{2^{i}j}+1\equiv 2\pmod 8$, whence
each of those $a-1$ factors contributes exactly $v_2=1$. Adding up,
\[
  v_2((m-1)^{n+2}-1)=v_2(m-2)+v_2(m)+a-1,
\]
and since $d=((m-1)^{n+2}-1)/m$ the formula follows. 
\end{proof}

\begin{thm}\label{thmAIV}
$\cE^n(\cK)$ is finite. More precisely, if $e=p_1^{\alpha_1}\cdots p_k^{\alpha_k}\in\cE^n_m(\cK)$, define for every $i=1,\ldots,k$
$$
G_i:=(\Z/p_i^{\alpha_i})^\times/\langle -1,1-m\rangle.
$$
Then 
$$
A:=\{1\le i\le k:|G_i|=1\}\subseteq\{1\le i\le k:\varphi(p_i^{\alpha_i})\mid n+2\},
$$
$$
B:=\{1\le j\le k:|G_j|\neq 1\}\subseteq\{1\le j\le k: p_j^{\alpha_j}\mid p_i^{(n+2)\rho}-1\text{ for some }i\in A\text{, and }\rho\le 2|A|\}.
$$
Furthermore, for every $j\in B$ we have the inequality
\begin{equation}\label{eq3.2}
1\le \frac{1}{|G_j|}+\sum_{i\in A}\frac{1}{\mathrm{ord}_{G_j}(p_i)}.
\end{equation}
\end{thm}

\begin{proof}
Note first that if $i\in A$, then $(\Z/p_i^{\alpha_i})^\times$ is generated by $-1$ and $1-m$. If $(\Z/p_i^{\alpha_i})^\times$ is cyclic, then it is generated by an element of the form $r=\pm(1-m)^h$. Hence $\varphi(p_i^{\alpha_i})=\text{ord}_{p_i^{\alpha_i}}(r)\mid \text{ord}_{e}(r)\mid n+2$. This proves the assertion about $A$ in this case. Otherwise, $p_i=2$ and $\alpha_i=\alpha\ge 3$, then $(\Z/2^{\alpha})^\times\simeq (\Z/2\Z)\oplus(\Z/2^{\alpha-2})$ is generated by $-1$ and $5$. Hence $m-1\equiv_{2^\alpha}\pm5^u$ for some odd $u$, and so $[m-1]_8\in\{3,5\}$ which implies that $v_2(m-2)\le 2$. Then, by \cref{lemma}
$$
\alpha=v_2(e)\le v_2(d)=v_2(m-2)+v_2(n+2)-1\le v_2(n+2)+1, 
$$
and so $\varphi(2^\alpha)=2^{\alpha-1}\mid n+2$. This completes the assertion about $A$. For the assertions about $B$ consider the set
$$
\Psi_e:=\{\chi\in C^-(e): \chi(1-m)=1\}.
$$
We claim that $\Psi_e\neq\varnothing$, in fact note that $-1\notin\langle1-m\rangle\subseteq(\Z/e)^\times$ since otherwise $e\mid (1-m)^h+1$ for some $h\mid \frac{n}{2}+1$ and so, by \cref{prop2.2}, $\frac{d}{e}\in\fD^n_m(\cK)$ contrary to the fact that $e\in \cE^n_m(\cK)$. Fix now some $\chi\in\Psi_e$ of maximal conductor. Let $K:=\langle -1,1-m\rangle\subseteq(\Z/e)^\times$ and $K_i:=\langle-1,1-m\rangle\subseteq(\Z/p_i^{\alpha_i})^\times$ for every $i=1,\ldots,k$. Fix some $j\in B$ and consider any character $\theta:(\Z/p_j^{\alpha_j})^\times\to\C^\times$ such that $K_j\subseteq\ker(\theta)$, then $\chi\theta\in\Psi_e$ and
$$
\text{cond}(\chi\theta)=\prod_p\text{cond}(\chi_p\theta_p)=\text{cond}(\chi_{p_j}\theta)\cdot\prod_{p\neq p_j}\text{cond}(\chi_p).
$$
If $\chi_{p_j}=1$, then by the maximality of the conductor we have $\text{cond}(\theta)=1$ and so $(\Z/p_j^{\alpha_j})^\times/K_j=G_j$ has no non-trivial characters, i.e. it is a trivial group, then $j\in A$ which is a contradiction. Hence $\chi_{p_j}\neq 1$. By the maximality of the conductor
$$
g:=\text{cond}(\chi\theta)\mid \text{cond}(\chi)=\prod_{\chi_{p_r}\neq 1}p_r^{\alpha_r}.
$$
Let $f:=e/g$, then $\chi\theta\in PC^-(e/f)$ with $(\chi\theta)(1-m)=1$. By \cref{prop3.1} there exists some $p_i\mid f=e/g$ and $p_i\nmid e/f=g$ such that $(\chi\theta)(p_i)=1$. Since the only possible difference between $g$ and $\text{cond}(\chi)$ is a power of $p_j$, this implies that $i=j$ or $\chi_{p_i}=1$, and so  $i\in A\cup\{j\}$ and if $i=j$ then $p_j\nmid g=\text{cond}(\chi\theta)$ and so $\chi_{p_j}\theta=1$, i.e. $\theta=\chi_{p_j}^{-1}$. On the other hand, if $i\in A$ then
$$
\theta(p_i)=\chi(p_i)^{-1}=:c_i
$$
has a fixed value at $p_i$. This shows that
$$
\widehat{G_j}\subseteq\bigcup_{i\in A}\{\theta\in G_j: \theta(p_i)=c_i\}\cup\{\chi_{p_j}^{-1}\}=\bigcup_{i\in A}C_i\cup\{\chi_{p_j}^{-1}\}
$$
where each $C_i$ has a transitive multiplicative action of $\widehat{G_j/\langle p_j\rangle}$. Hence 
$$
|G_j|\le 1+\sum_{i\in A}\frac{|G_j|}{\text{ord}_{G_j}(p_i)}
$$
and so \eqref{eq3.2} follows. Finally, since $|G_j|\ge 2$, we see that in order to satisfy \eqref{eq3.2} there must exists some $i\in A$ such that $\rho:=\text{ord}_{G_j}(p_i)\le 2|A|$, then $p_i^\rho\equiv_{p_j^{\alpha_j}}\pm(1-m)^h$ and so $p_i^{(n+1)\rho}\equiv_{p_j^{\alpha_j}}1$. This proves the assertions about $B$. From both bounds on $A$ and $B$, which are independent of $m$, we have bound the possible prime power factors of $e$ and so we have the finiteness of $\cE^n(\cK)$.
\end{proof}

The previous theorem together with \cref{rmk3.1} imply \cref{thm0.2}. Let us now prove \cref{thm0.3}.

\bigskip

\begin{proof}[of \cref{thm0.3}]
Assume that for all primes $p$ such that $p-1\mid n+2$ it holds that $p\nmid d$. Suppose that $\cE^n_m(\cK)\neq\varnothing$, take $e\in\cE^n_m(\cK)$. In the notation of \cref{thmAIV}, if $A=\varnothing$, then $B=\varnothing$ and so $e$ has no prime factors, which is absurd. Then $A\neq\varnothing$, and so there exists some $i\in A$ such that $p_i\mid e$ and $\varphi(p_i^{\alpha_i})=p_i^{\alpha_i-1}(p_i-1)\mid n+2$. Then $p_i-1\mid n+2$ and $p_i\mid d$ which is a contradiction to our hypothesis. Therefore $\cE^n_m(\cK)=\varnothing$ and so $\fB^n_m(\cK)=\fD^n_m(\cK)$.
\end{proof}

\section{Hodge cycles on Klein varieties of small dimension}
\label{sec4}
In this final section we describe the set of exceptional denominators for small dimensions. By \eqref{eq3.1} the determination of $\fB_m^n(\cK)$ for all $m$
amounts to computing the set $\cE^n(\cK)$ together with the residue sets $S_e$. Iterating over $e$, by \cref{rmk3.1} it
suffices, for each $e$, to run over one complete residue system for $m$ modulo
$e$ and test
\begin{equation}\label{eq:search}
e\mid d,\qquad e\nmid (1-m)^h+1\ \ \text{for all }h\text{ admissible},\qquad
\alpha_0\in\fB_e^n,
\end{equation}
the last condition being a finite check on a single
vector modulo $e$, and $h$ is called admissible if $2h\mid n+2$ and $n+2=2^th$. This yields $\cE_m^n(\cK)$ and the sets $S_e$ for all $e$ up to
the search bound.

\begin{rmk}
In \cref{rmk3.1} was required that $e\nmid (1-m)^h+1$ for all $h\mid \frac{n}{2}+1$, but in \eqref{eq:search} we reduced the test only to $h$ admissible. This test reduction can be done because if $(1-m)^h\equiv_e-1$ for some $h\mid \frac{n}{2}+1$, then taking $\ell$ odd such that $n+2=2^th\ell$ we have that $h\ell$ is admissible and $(1-m)^{h\ell}\equiv_e-1$.
\end{rmk}

\begin{thm}\label{thm:computations}
A search using \eqref{eq:search} over $e\leq 10^7$ for $n=2,4,6$ and over $e\le 10^6$ for $n=8$ yields
\[
\renewcommand{\arraystretch}{1.5}
\begin{array}{rcl}
\cE_m^4(\cK)&\supseteq&\Big\{14_{([m]_{14}=4,6)},\ 21_{([m]_{21}=3,12)},\
 28_{([m]_{28}=6,18)},\ 39_{([m]_{39}=5,11)},\\
&&\ \ 42_{([m]_{42}=6,12,18,20,24,32)},\
 78_{([m]_{78}=18,24,30,36,44,50)},\ 84_{([m]_{84}=6,18)},\\
&&\ \ 126_{([m]_{126}=32,62)},\ 156_{([m]_{156}=18,102)},\
 168_{([m]_{168}=18,90)},\ 210_{([m]_{210}=102,132)}\Big\},\\
\cE_m^6(\cK)&\supseteq&\Big\{20_{([m]_{20}=4,8)},\ 30_{([m]_{30}=8,14)},\
 60_{([m]_{60}=8,14,38,44)},\ 102_{([m]_{102}=20,26,44,50)},\\
&&\ \ 120_{([m]_{120}=14,38)},\ 170_{([m]_{170}=54,78,84,128)},\
 204_{([m]_{204}=20,44,128,152)},\\
&&\ \ 510_{([m]_{510}=128,224,248,254)}\Big\},\\
\cE_m^8(\cK)&\supseteq&\Big\{33_{([m]_{33}=8,14,20,29)},\ 55_{([m]_{55}=5,10,15,50)},\
 62_{([m]_{62}=16,24,28,30)},\\
&&\ \ 66_{([m]_{66}=8,14,20,62)},\
 110_{([m]_{110}=10,20,30,40,42,50,52,60,62,70,80,102)},\\
&&\ \ 124_{([m]_{124}=30,78,86,90)},\ 186_{([m]_{186}=86,92,140,152)},\
 220_{([m]_{220}=30,130,150,190)},\\
&&\ \ 330_{([m]_{330}=20,80,140,260)},\ 363_{([m]_{363}=41,95,113,119)},\
 682_{([m]_{682}=214,334,340,512)},\\
&&\ \ 726_{([m]_{726}=404,458,476,482)},\
 2046_{([m]_{2046}=512,896,1016,1022)}\Big\}.
\end{array}
\]
\end{thm}

\begin{proof}
All computations reported were carried out with a short program
implementing \eqref{eq:search} directly. The only nontrivial ingredient is the
test $\alpha_0\in\fB_e^n$, which requires evaluating
$|t\cdot\alpha_0|$ for the $\varphi(e)$ units $t$ modulo $e$.
\end{proof}

\begin{rmk}
The same search over $e\le 10^6$ in dimension $n=10$ gives
\[
\begin{array}{rcl}
\cE^{10}(\cK)&\supseteq&\{14,20,21,28,30,39,42,52,60,70,78,84,90,105,120,126,140,155,
156,168,180,\\
&&182,195,210,234,280,310,360,390,420,546,620,630,910,1170,1638,2730,8190\},
\end{array}
\]
together with the corresponding residue sets $S_e$. We are led to the following expectation.
\end{rmk}

\begin{conj}\label{conj:complete}
The inclusions of \cref{thm:computations} are equalities for $n=4,6,8$ and
any $m$.
\end{conj}

The inclusion of \cref{thm:computations} is an equality for $n=2$ and any $m$, since by \cite[Theorem C]{Aoki1983} the exceptional divisors in $\cE^2(\cK)$ are bounded by $180$. For $n\geq4$ the bounds obtained by the proof of \cref{thm0.2} are astronomical, and so they are useless for computer implementation. Nevertheless, the fact that only $5$ of the $22080$
exceptional characters of Fermat surfaces occur as test vectors
$\alpha_0\in\fB^2_e$ suggests that the constraint imposed by the Klein covering is
severe, and one may hope that it forces $\cE^n(\cK)$ to be smaller than $10^4$ for small dimensions $n\le 10$. We record a compatibility which explains
several coincidences in the data above. Write $d^{(n)}(m)=\big((m-1)^{n+2}-1\big)/m$
and $\alpha_0^{(n)}=\big(1,1-m,\dots,(1-m)^{n+1}\big)$ for the invariants of
dimension $n$.

\begin{prop}\label{prop:transfer}
Let $n,n'$ be even with $n+2\mid n'+2$, and let $e\geq2$ divide
$d^{(n)}(m)$. Then $e$ divides $d^{(n')}(m)$, and
\[
\alpha_0^{(n)}\in\fB_e^{n}
\quad\Longrightarrow\quad
\alpha_0^{(n')}\in\fB_e^{n'}.
\]
Consequently, if $e\in\cE_m^n(\cK)$ and $e\nmid (1-m)^h+1$ for all $h$ admissible,
then $e\in\cE_m^{n'}(\cK)$.
\end{prop}

\begin{proof}
Put $q=(n'+2)/(n+2)$. Since $e\mid d^{(n)}(m)\mid (m-1)^{n+2}-1$ we have
$(m-1)^{\,n+2}\equiv_e1$, so $(m-1)^{\,n'+2}\equiv_e1$ and
$(m-1)^{\,n'+2}-1=((m-1)^{\,n+2}-1)\sum_{t<q}(m-1)^{\,t(n+2)}$ gives
$d^{(n)}(m)\mid d^{(n')}(m)$. The same congruence shows that
$\alpha_0^{(n')}$ reduces modulo $e$ to the $q$-fold concatenation of
$\alpha_0^{(n)}$, whence for every $t\in(\Z/e\Z)^\times$
\[
\big|t\cdot\alpha_0^{(n')}\big|
=q\,\big|t\cdot\alpha_0^{(n)}\big|
=q\Big(\frac n2+1\Big)=\frac{n'}2+1 .
\]
The last assertion
follows directly from this.
\end{proof}

\begin{rmk}
The data are consistent with the previous proposition in every instance
where it applies. Since $4\mid8$ one predicts $\cE^2(\cK)\subseteq\cE^6(\cK)$, and indeed
$\{20,30,60,120\}\subseteq\cE^6(\cK)$. Since $4\mid12$ and
$6\mid12$ one predicts $\cE^2(\cK)\cup\cE^4(\cK)\subseteq\cE^{10}(\cK)$, and all fifteen
denominators do occur. Conversely $4\nmid6$, $8\nmid12$ and
$10\nmid12$, and correspondingly $\cE^2(\cK)\not\subseteq\cE^4(\cK)$, while none of
$102,170,204,510$ or of the dimension-$8$ denominators appears in $\cE^{10}(\cK)$. 
\end{rmk}

Finally we prove the last result of the article.

\bigskip

\begin{proof}[of \cref{cor0.5}]
By \cref{rmk2.1}, the Hodge conjecture for $Y^n_m$ is reduced to show that $\alpha_0\in \fC^n_{e}$ for all $e\in\cE^n_m(\cK)$. In particular, since the Hodge conjecture holds for Fermat surfaces and fourfolds \cite{kang2016refined}, the result follows for $n=2,4$. On the other hand, for $n=6$ and $m\le 10^7$, by \cref{thm0.4} we have to verify that $\alpha_0$ is algebraic for each $e\in \cE^6_m(\cK)$. By \cite[Theorem 1.1]{MMRV2026} we cover all $e$ except by $e\in\{102, 170, 204, 510\}$. On the other hand, applying the length reduction algorithm of \cite[Algorithm 1]{MMRV2026} we obtain the following decomposition of $\alpha_0$ in terms of linear, standard and exceptional tuples:

\begin{table}[ht]
\centering
\begin{tabular}{@{}ll@{}ll@{}}
\toprule
$e$ & $S_e$ &\ \ \ \textbf{Decomposition of }$\alpha_0$ \\
\midrule
$102=2\cdot 3\cdot 17$ & 20, 26, 44, 50 &\ \ \ \ $D+s_2+s_3+s_{17}$ \\
\midrule
$170=2\cdot 5\cdot 17$ & 54, 78, 84, 128 &\ \ \ \ $D+s_2+s_5+s_{17}$ \\
\midrule
$204=2^2\cdot 3\cdot 17$ & 20, 44, 128, 152 &\ \ \ \ $D+s_2+s_3+\xi_{68}$ \\
\midrule
$510=2\cdot3\cdot5\cdot 17$ & 128, 224, 248, 254 &\ \ \ \  $D+s_2+s_3+s_5+s_{17}$ \\
\bottomrule
\end{tabular}
\end{table}

We remark that for every $e$, the tuples $\alpha_0$ obtained for all admissible values of $[m]_e\in S_e$ are the same up to a permutation (and so have the same decomposition). In particular, except by $e=204$, we have that $\alpha_0\in\fC^6_e$. Doing a similar analysis for $n=8$ and $m\le 10^6$, by \cite[Theorem 1.1]{MMRV2026} we discard all $e$ except $e\in\{124, 186,220,682,2046\}$. Running \cite[Algorithm 1]{MMRV2026} for them (except $e=2046$ due to our computer limitations) we obtain the following decompositions:

\begin{table}[ht]
\centering
\begin{tabular}{@{}ll@{}ll@{}}
\toprule
$e$ & $S_e$ &\ \ \ \textbf{Decomposition of }$\alpha_0$ \\
\midrule
$124=2^2\cdot 31$ & 30, 78, 86, 90 &\ \ \ \ $D+s_2$ \\
\midrule
$186=2\cdot 3\cdot 31$ & 86, 92, 140, 152 &\ \ \ \ $D+s_2+s_3+s_{31}$ \\
\midrule
$220=2^2\cdot 5\cdot 11$ & 30, 130, 150, 190 &\ \ \ \ $D+s_2+s_5+s_{11}+\xi_{55}$ \\
\midrule
$682=2\cdot11\cdot31$ & 214, 334, 340, 512 &\ \ \ \  $D+s_2+s_{11}+s_{31}$ \\
\bottomrule
\end{tabular}
\end{table}

In consequence, we have proved that $\alpha_0\in\fC^8_e$ for all $e$ except by $e\in\{220,2046\}$.
\end{proof}

\bibliographystyle{alpha}
\bibliography{biblio}

\end{document}